\documentclass[12pt,twoside]{amsart}

\usepackage{amsmath}
\usepackage{amssymb}
\usepackage{amscd}

\DeclareMathSymbol{\ordcol}{\mathord}{operators}{'072}
\newcommand{\col}{{{\hskip1.5pt\ordcol\hskip1.5pt}}}

\renewcommand{\bar}{\overline}

\newcommand{\CC}{\mathbb{C}}

\newcommand{\FF}{\mathbb{F}}

\newcommand{\PP}{\mathbb{P}}
\newcommand{\QQ}{\mathbb{Q}}
\newcommand{\RR}{\mathbb{R}}

\newcommand{\Cv}{\CC_v}

\newcommand{\ints}{{\mathcal O}}

\newcommand{\maxid}{{\mathcal M}}

\newcommand{\calZ}{{\mathcal Z}}

\newcommand{\PCv}{\PP^1(\Cv)}

\newcommand{\Pkhat}{\PP^1(\hat{k})}
\newcommand{\PKhat}{\PP^1(\hat{K})}
\newcommand{\PK}{\PP^1(K)}

\newcommand{\hhat}{\hat{h}}
\newcommand{\Khat}{\hat{K}}

\newcommand{\frakK}{{\mathfrak{K}}}

\newcommand{\PGL}{\hbox{\rm PGL}}

\newcommand{\Dbar}{\bar{D}}

\theoremstyle{plain}
\newtheorem{thm}{Theorem}[section]
\newtheorem{prop}[thm]{Proposition}

\newtheorem{lemma}[thm]{Lemma}
\newtheorem{defin}[thm]{Definition}

\newtheorem*{thmA}{Theorem A}
\newtheorem*{thmB}{Theorem B}

\theoremstyle{definition}

\newtheorem{remark}[thm]{Remark}

\title[Preperiodic points over function fields]
      {Heights and preperiodic points of polynomials over function fields}
\author{Robert L. Benedetto}
\date{October 20, 2005; revised December 12, 2005}
\thanks{The author gratefully acknowledges the support
  of a Miner D.\ Crary Research Fellowship from Amherst College
  and NSA Young Investigator Grant H98230-05-1-0057}
\subjclass[2000]{Primary: 11G50.  Secondary: 11D45, 37F10}
\keywords{filled Julia set, canonical height}
\address{Department of Mathematics and Computer Science \\
        Amherst College \\
        Amherst, MA 01002 \\
        USA}
\email{rlb@cs.amherst.edu}
\urladdr{http://www.cs.amherst.edu/\textasciitilde rlb}

\begin{document}

\newcounter{bean}
\newcounter{sheep}

\begin{abstract}
Let $K$ be a function field in one variable over an arbitrary field
$\FF$.  Given a rational function $\phi\in K(z)$ of degree at least
two, the associated canonical height on the projective line was
defined by Call and Silverman.  The preperiodic points of $\phi$ all
have canonical height zero; conversely, if $\FF$ is a finite field,
then every point of canonical height zero is preperiodic.  However, if
$\FF$ is an infinite field, then there may be non-preperiodic points
of canonical height zero.  In this paper, we show that for polynomial
$\phi$, such points exist only if $\phi$ is isotrivial.
In fact, such $K$-rational points exist only if $\phi$ is
defined over the constant field of $K$
after a $K$-rational change of coordinates.
\end{abstract}

\maketitle

Let $K$ be a field with algebraic closure $\Khat$,
and let $\phi:\PKhat\rightarrow\PKhat$ be
a morphism defined over $K$.  We may write $\phi$ as a
rational function $\phi\in K(z)$.  Denote the $n^{\text{th}}$
iterate of $\phi$ under composition by $\phi^n$.  That is,
$\phi^0$ is the identity function, and for $n\geq 1$,
$\phi^n = \phi\circ \phi^{n-1}$.  A point $x\in\PKhat$
is said to be
{\em preperiodic} under $\phi$ if there are integers $n>m\geq 0$
such that $\phi^m(x)=\phi^n(x)$.
Note that $x$ is preperiodic if and only
if its forward orbit $\{\phi^n(x):n\geq 0\}$ is finite.

If $K$ is a number field or a function field in one variable,
and if $\deg\phi \geq 2$,
there is a {\em canonical height} function
$\hhat_{\phi}:\PKhat\rightarrow \RR$ associated to $\phi$.
(In this context, the degree $\deg\phi$ of $\phi$
is the maximum of the the degrees of its numerator and denominator.)
The canonical height gives a rough measure of the arithmetic
complexity of a given point, and it also satisfies the functional
equation $\hhat_{\phi}(\phi(x)) = d\cdot \hhat_{\phi}(x)$, where
$d=\deg\phi$.
Canonical heights will be discussed further in Section~\ref{sect:heights};
for more details, we refer the reader to the original paper
\cite{CS1} of Call and Silverman, or to the exposition in
Part~B of \cite{HSbook}.

All preperiodic points of $\phi$ clearly have canonical height zero.
Conversely, if $K$ is a global field (i.e., a
number field or a function field in one variable over a finite field),
then $\hhat_{\phi}(x)=0$
if and only if $x$ is preperiodic.  This equivalence is
very useful in the study of rational preperiodic points
over such fields; see, for example, \cite{Ben14,CaGol}.

However, if $K$ is a function field over
an infinite field, there may be points of canonical height
zero which are not preperiodic.  For example, if $K=\QQ(T)$
and $\phi(z)=z^2$, then the preperiodic points in $\PKhat$
are $0$, $\infty$, and the roots of unity in $\hat{\QQ}$; however,
{\em all} points in $\PP^1(\hat{\QQ})$ have canonical height
zero.  Similarly, for the same field $K$, consider the function
$\psi(z) = T z^3$.  In this case, the only $K$-rational points
of canonical height zero are $0$ and $\infty$, both of which are
preperiodic.  Nevertheless, in
$\PKhat$, any point of the form $a T^{-1/2}$ with $a\in\hat{\QQ}$
has canonical height zero, but such a point is preperiodic if and
only if $a$ is either $0$, $\infty,$ or a root of unity.

It is important to note that both of the examples in the
previous paragraph are of
isotrivial maps (see Section~\ref{sect:jul}).
Specifically, $K=\QQ(T)$ has constant field $\QQ$, and
$\phi$ is defined over $\QQ$ as written.  Similarly, although
the coefficients of $\psi$ are not constants of $K$,
the $K(T^{1/2})$-rational
change of coordinates $\gamma(z) = T^{-1/2} z$ makes
$\psi$ into the map
$\gamma^{-1}\circ\psi\circ\gamma(z)=z^3$, which
is defined over $\QQ$.

In the theory of elliptic curves over function fields,
it is well known that certain finiteness results for heights
hold when the curve is not isotrivial;
see, for example, Theorem~5.4 of \cite{SilAT}.
In the same vein,
the main result of this paper is that in the dynamical setting,
once isotrivial maps
are excluded, canonical heights on one-variable function fields
have the same relation to preperiodic points as their analogues
on global fields.

\begin{thmA}
Let $\FF$ be an arbitrary field, let $K$ be a finite extension
of $\FF(T)$, and let 
$\phi\in K[z]$ be a polynomial of degree at least two.
Suppose that there is no $K$-rational
affine change of coordinates $\gamma(z)=az+b$ for which
$\gamma^{-1}\circ\phi\circ\gamma$ is defined over the
constant field of $K$.
Then for every point $x\in\PK$, $x$ is preperiodic
under $\phi$ if and only if $\hhat_{\phi}(x)=0$.
\end{thmA}

The substance of Theorem~A is the statement
that set $\calZ$ of $K$-rational
points of canonical height zero is finite.
In fact, as we will observe in Remark~\ref{rem:slogs},
under the hypotheses of Theorem~A,
we can bound the size of $\calZ$
in terms of the number $s$ of places of bad reduction.
Specifically, by invoking the arguments in \cite{Ben14},
we can obtain a bound of $O(s\log s)$.  Such statements
are related to Northcott's finiteness theorem \cite{Nor}
and Morton and Silverman's uniform boundedness conjecture
\cite{MorSil1}, both originally stated for global fields.

In addition, as an immediate consequence of Theorem~A, we have the
following result.

\begin{thmB}
Let $\FF$ and $K$ be as in Theorem~A, and let
$\phi\in K[z]$ be a polynomial of degree at least two.
Suppose that $\phi$ is not isotrivial.
Then for every point $x\in\PKhat$, $x$ is preperiodic
under $\phi$ if and only if $\hhat_{\phi}(x)=0$.
\end{thmB}

To prove the theorems, we will recall some basic background
on function fields and their associated
local fields in Section~\ref{sect:nota}.  Our main tools
will be the standard notions of good and bad reduction
and of filled Julia sets, defined in Section~\ref{sect:jul}.
In Section~\ref{sect:heights}, we will recall the
theory of canonical heights,
which we will relate to reduction and filled Julia sets.
Finally, after stating some lemmas in Section~\ref{sect:lem}, 
we will prove our theorems in Section~\ref{sect:pfs}.

The author would like to thank Lucien Szpiro and Thomas Tucker
for introducing him to the questions considered in this paper.
Many thanks to Matthew Baker, Joseph Silverman, and again to
Thomas Tucker, as
well as to the referee, for a number of helpful suggestions
to improve the exposition.



\section{Background on Function Fields and Local Fields}
\label{sect:nota}

In this section we recall the necessary fundamentals
from the theory of function fields and their associated local fields.
We also set some notational conventions for this paper.
For more details on function fields and their absolute
values, we refer the reader to 
Chapter~5 of \cite{Rosen}, as well as
Chapters~1--2 of \cite{Lang2},
Section~B.1 of \cite{HSbook}, and Section~4.4 of \cite{RV}.
See \cite{Gou,Kob} for expositions concerning the local
fields $\Cv$.

\subsection{Function Fields, Places, and Absolute Values}
\label{ssect:fields}
We set the following notation throughout this paper.
\begin{tabbing}
\hspace{0.3in} \= $\FF$ \= \hspace{0.3in} \= an arbitrary field \\
\> $K$ \> \> a function field in one variable over $\FF$;
i.e., a finite extension of $\FF(T)$ \\
\> $\hat{K}$ \> \> an algebraic closure of $K$ \\
\> $M_K$ \> \> a proper set of absolute values on $K$, satisfying
            a product formula.
\end{tabbing}
We will explain this terminology momentarily.
Here, $\FF(T)$ is the field of rational functions in one variable
with coefficients in $\FF$.
Recall that the {\em constant field} $\FF_K$ of $K$ is the set
of all elements of $K$ that are algebraic over $\FF$.  It is a finite
extension of $\FF$.

Recall that an {\em absolute value} on $K$ is a function
$|\cdot|_v:K\rightarrow \RR$ satisfying
$|x|_v\geq 0$ (with equality if and only if $x=0$),
$|xy|_v = |x|_v |y|_v$, and
$|x+y|_v\leq |x|_v + |y|_v$, for all $x,y\in K$.
The {\em trivial} absolute value is given by $|0|_v=0$
and $|x|_v=1$ for $x\in K^{\times}$.
Two absolute values on $K$ are said to be {\em equivalent} if
they induce the same topology on $K$.

In this setting,
to say that the set $M_K$ is {\em proper} is to say that each $v\in M_K$
is nontrivial; that $v,w\in M_K$ are equivalent if and only if $v=w$;
and that for any $x\in K^{\times}$, there are only finitely
many $v\in M_K$ for which $|x|_v\neq 1$.
That $M_K$ satisfies a {\em product formula} is to say that
there are positive real numbers $\{n_v:v\in M_K\}$ such that for
any $x\in K^{\times}$,
\begin{equation}
\label{eq:pfmla}
\prod_{v\in M_K} |x|_v^{n_v} = 1.
\end{equation}
(Note that the properness of $M_K$ means that for any given $x$,
the product on the left-hand side of \eqref{eq:pfmla} is really
a finite product.)
The absolute values $v\in M_K$ are called the {\em places} of $K$,
though we note that some authors consider a place to be
an equivalence class of absolute values on $K$.

The existence of the set $M_K$ in the function field setting
is guaranteed by the association between primes
and equivalence classes of nontrivial absolute values on $K$.
(Here, a prime is the maximal ideal $P$
of a discrete valuation ring $R\subseteq K$
such that $\FF_K\subseteq R$ and
$K$ is the field of fractions of $R$.)

Because $K$ is a function field, we also have the following
additional facts.  First,
all absolute values $v$ in $M_K$
are {\em non-archimedean}, which is to say that they satisfy
the ultrametric triangle inequality
$$|x+y|_v \leq \max\{|x|_v,|y|_v\}.$$
Second, for any $x\in K^{\times}$, we have
$x\in \FF^{\times}_K$ (i.e., $x$ is a nonzero constant) if and only if
$|x|_v=1$ for all $v\in M_K$.
Third, we may
normalize our absolute values so that $n_v=1$ for all $v\in M_K$
in \eqref{eq:pfmla}, giving the simplified product formula
$$\prod_{v\in M_K} |x|_v = 1
\qquad \text{ for all } x\in K^{\times}
$$

\subsection{Local fields}
For each $v\in M_K$, we can form the local field $K_v$,
which is the completion of $K$ with respect to $|\cdot|_v$.
We write $\Cv$ for the completion of an algebraic
closure $\Khat_v$ of $K_v$.  (The absolute value $v$ extends in
a unique way to $\Khat_v$ and hence to $\Cv$.)
The field $\Cv$ is then a complete and algebraically closed field.

Because $v$ is non-archimedean, the disk
$\ints_v = \{c\in\Cv : |c|_v \leq 1\}$ forms a ring,
called the {\em ring of integers} of $\Cv$, which has a unique
maximal ideal
$\maxid_v = \{c\in\Cv : |c|_v < 1\}$.  The quotient
$k_v=\ints_v/\maxid_v$ is called the {\em residue field} of $\Cv$.
The natural reduction map from $\ints_v$ to $k_v$, sending
$a\in\ints$ to $\bar{a}=a+\maxid_v\in k_v$, will
be used to define good and bad reduction of a polynomial
in Definition~\ref{def:reduc} below; but after invoking a few
simple lemmas about good and bad reduction, we will
not need to refer to $\ints_v$, $\maxid_v$, or $k_v$ again.

\subsection{Disks}
Let $\Cv$ be a complete and algebraically closed field with
absolute value $v$.  Given $a\in\Cv$ and $r>0$, we write
$$\Dbar(a,r) = \{x\in\Cv : |x-a|_v \leq r\}
\quad\text{and}\quad
D(a,r) = \{x\in\Cv : |x-a|_v < r\}$$
for the closed and open disks, respectively, of radius $r$
centered at $a$.  Note our convention that all
disks have positive radius.

If $v$ is non-archimedean and $U\subseteq\Cv$ is a disk,
then the radius of $U$ is unique; it is the same as the diameter
of the set $U$ viewed as a metric space.  However, any point
$b\in U$ is a center.  That is, if $|b-a|_v\leq r$, then
$\Dbar(a,r)=\Dbar(b,r)$, and similarly for open disks.
It follows that two disks intersect if and only
if one contains the other.


\section{Bad Reduction and Filled Julia Sets}
\label{sect:jul}

The following definition originally appeared in \cite{MorSil1}.
We have modified it slightly so that ``bad reduction'' now
means not potentially good, as opposed to not good.
\begin{defin}
\label{def:reduc}
Let $\Cv$ be a complete, algebraically closed
non-archimedean field with
absolute value $|\cdot|_v$, ring of integers
$\ints_v = \{c\in\Cv : |c|_v \leq 1\}$,
and residue field $k_v$.
Let $\phi(z) \in \Cv(z)$ be a rational function
with homogenous presentation
$$\phi\left( [x,y]\right) = [f(x,y),g(x,y)],$$
where $f,g\in\ints_v[x,y]$ are relatively prime homogeneous polynomials
of degree $d=\deg\phi$, and at least one coefficient of $f$ or $g$
has absolute value $1$.  We say that $\phi$ has {\em good reduction}
at $v$ if $\bar{f}$ and $\bar{g}$ have no common zeros in
$k_v \times k_v$ besides $(x,y)=(0,0)$.
We say that $\phi$ has {\em potentially good reduction} at $v$
if there is some linear fractional transformation
$h\in\PGL(2,\Cv)$ such that $h^{-1}\circ\phi\circ h$ has good
reduction.
If $\phi$ does
not have potentially good reduction,
we say it has {\em bad reduction} at $v$.
\end{defin}
Naturally, for $f(x,y)=\sum_{i=0}^d a_i x^i y^{d-i}$,
the reduction $\bar{f}(x,y)$ in Definition~\ref{def:reduc}
means $\sum_{i=0}^d \bar{a}_i x^i y^{d-i}$.

In this paper, we will consider only polynomial $\phi$
of degree at least $2$; that is, $\phi(z)=a_d z^d + \cdots + a_0$,
where $d\geq 2$, $a_i\in\Cv$, and $a_d\neq 0$.
It is easy to check that such $\phi$ has good reduction
if and only if $|a_i|_v\leq 1$ for all $i$ and $|a_d|_v = 1$.
In particular, by the properness of $M_K$,
there are only finitely many places $v\in M_K$ at which
$\phi$ has bad reduction.
(The same finiteness result also holds for $\phi\in K(z)$,
but we do not need that fact here.)

We recall the following standard definition from polynomial dynamics.
\begin{defin}
\label{def:julia}
Let $\Cv$ be a complete, algebraically closed
field with absolute value $|\cdot|_v$, and let
$\phi(z) \in \Cv[z]$ be a polynomial of degree $d\geq 2$.
The {\em filled Julia set} of $\phi$ at $v$ is
$$\frakK_{\phi,v}
= \{x\in\Cv: \{|\phi^n(x)|_v\}_{n\geq 0} \text{ is bounded} \}.$$
\end{defin}

We recall the following three fundamental facts.
First, $\frakK_{\phi,v}$ is invariant
under $\phi$; that is,
$\phi^{-1}(\frakK_{\phi,v}) = \phi(\frakK_{\phi,v})=\frakK_{\phi,v}$.
Second, all the finite preperiodic points of $\phi$
(that is, all the preperiodic points in $\PCv$ other than the
fixed point at $\infty$) are contained in $\frakK_{\phi,v}$.
Finally, the polynomial $\phi\in\Cv[z]$ has good reduction
if and only if $\frakK_{\phi,v}=\Dbar(0,1)$.
(However, if $\phi$ has bad reduction, then $\frakK_{\phi,v}$
can be a very complicated fractal set.)

Filled Julia sets have been studied extensively
in the archimedean case $\Cv=\CC$; see, for example,
\cite{Bea,CG,Mil}.  For the non-archimedean setting,
see \cite{Ben5,Ben14,Riv1}.

Finally, we recall that if $\phi\in K(z)$ is a rational function
over the function field $K$, we say that $\phi$ is {\em isotrivial}
if there is a finite extension $L/K$ and an $L$-rational
change of coordinates $\gamma\in\PGL(2,L)$ such that
$\gamma^{-1}\circ \phi \circ \gamma$ is defined over
a the constant field $\FF_L$ of $L$.
Note that if $\phi$ is a polynomial and such a $\gamma$ exists,
then there is a map $\beta$ of the form $\beta(z)=az + b$
with $a\in L^{\times}$ and $b\in L$ such that
$\beta^{-1}\circ \phi \circ \beta$ is defined over $\FF_L$.

We comment for the interested reader that the term ``isotrivial''
is borrowed from the theory of moduli spaces.  In that context,
a trivial family has no dependence on the parameter(s), and
an isotrivial family is one where all smooth fibers are isomorphic.
See \cite{Beau}, for example, for more information.

\section{Canonical Heights}
\label{sect:heights}

We recall the following facts from the theory of
heights and
canonical heights from \cite{CS1,HSbook,Lang2}.
The standard {\em height function} $h:\PKhat\rightarrow [0,\infty)$
is given by the formula
\begin{equation}
\label{eq:htdef1}
h(x) = \frac{1}{[L:K]}
\sum_{v\in M_K} \sum_{w\in M_L, \, w \mid v}
n_w \log \max\{|x_0|_w, |x_1|_w\},
\end{equation}
where $L/K$ is a finite extension over which $x=[x_0,x_1]$
is defined, $n_w=[L_w:K_v]$ is the corresponding local field
extension degree at $w$,
and $w\mid v$ means that
the restriction of $w$ to $K$ is $v$.
It is well known that \eqref{eq:htdef1}
is independent of the choices of homogeneous coordinates $x_0,x_1$
for $x$ and of the extension $L$.  In the special case that
$x=f/g\in \FF(T)$,
where $f,g\in\FF[T]$ are relatively prime polynomials, \eqref{eq:htdef1}
can be rewritten as
$$
h(f/g) = \max\{ \deg f , \deg g\} ,
$$
assuming a certain standard normalization is chosen for the
absolute values.

In addition, if $\FF$ is finite, then $h$ is {\em nondegenerate}
in the sense that for any real constants $B,D$, the set of $x\in\PKhat$
with $h(x)\leq B$ and $[K(x)\col K \hspace{0.75pt}]\leq D$
is a finite set.  Conversely, if $\FF$ is infinite, then $h$ fails
to be nondegenerate, because all (infinitely many) points of
$\PP^1(\FF)$ have height zero.  The essential
point of this paper is that certain finiteness results can still be proven
under weak hypotheses even in the case that $\FF$ is infinite.

One crucial property of the height function is that it
satisfies an approximate functional equation for any
morphism.  Indeed,
if $\phi\in K(z)$ is a rational function of degree $d\geq 1$,
there is a constant $C=C_{\phi}\geq 0$ such that
$$ \text{for all } x\in \PKhat,
\quad \left| h\left(\phi(x) \right) - d \cdot h(x)\right| \leq C .
$$

For a {\em fixed} $\phi(z)$ of degree $d\geq 2$,
Call and Silverman \cite{CS1} introduced a related {\em canonical}
height function
$\hat{h}_{\phi} : \Pkhat \rightarrow [0,\infty)$,
with the property that there
is a constant $C'=C'_{\phi}\geq 0$ such that
for all $x\in\PKhat$,
\begin{equation}
\label{eq:canht}
\hat{h}_{\phi}(\phi(x) ) = d \cdot \hat{h}(x)
\quad \text{and} \quad
| \hat{h}_{\phi}(x) - h(x) | \leq C' .
\end{equation}
(The height $\hhat_{\phi}$ generalized a construction
of N\'eron and Tate; see \cite{LangTate,Ner}.)
Note that by \eqref{eq:canht}, a preperiodic
point $x$ of $\phi$ must have canonical height zero.

Call and Silverman proved that
the canonical height has a decomposition
\begin{equation}
\label{eq:lchtsum}
\hhat_{\phi}(x) = 
\sum_{v\in M_K} \hat{\lambda}_{\phi,v}(x)
\qquad \text{for all } x\in K=\PK\setminus\{\infty\},
\end{equation}
where $\hat{\lambda}_{\phi,v}$ is the {\em local canonical height}
for $\phi$ at $v$.
For all places $v$ of good reduction, we have simply
$\hat{\lambda}_{\phi,v}(x) = \log\max\{|x|_v,1\}$, as
in the corresponding term of \eqref{eq:htdef1}.
More generally, if $\phi$ is a polynomial (still of degree $d\geq 2$),
we have
\begin{equation}
\label{eq:limdef}
\hat{\lambda}_{\phi,v}(x) =
\lim_{n\rightarrow\infty} d^{-n} \log\max\{|\phi^n(x)|_v, 1\}.
\end{equation}
It is easy to show that the limit converges and is nonnegative.

We will be able to avoid direct use of local canonical
heights, because of the crucial fact that if $\phi$ is a polynomial,
then $\hat{\lambda}_{\phi,v}(x)\geq 0$ for all $x\in K$, with
equality if and only if $x\in\frakK_{\phi,v}$.
(In fact, $\hat{\lambda}_{\phi,v}$ satisfies the functional
equation $\hat{\lambda}_{\phi,v}(\phi(x))=
d \cdot \hat{\lambda}_{\phi,v}(x)$ in the polynomial case.)
Thus, a point $x\in K \subseteq \PKhat$
has (global) canonical height zero if and only if 
$x\in\frakK_{\phi,v}$ for all $v\in M_K$.

For more information on local canonical heights, we refer
the reader to \cite{CS1} and to simpler results for polynomials
in \cite{CaGol}.
(The results of \cite{CaGol} assume that $K$ is a number field,
but all of the relevant proofs
go through verbatim in the function field case.)
In particular, the above description of $\hat{\lambda}_{\phi,v}$
comes from the algorithm in Theorem~5.3 of \cite{CS1}, a
simplified form of which appears as Theorem~3.1 of \cite{CaGol}.
Indeed, the description \eqref{eq:limdef} for polynomials is the same as
that given in Theorem~4.2 of \cite{CaGol},
and the connections
between local canonical heights and filled Julia sets are given
in Corollary~5.3 of the same paper.

\begin{remark}
It is easy to show, using the defining properties \eqref{eq:canht}
for $\hhat_{\phi}$, that a point $x\in\PKhat$ has canonical height
zero if and only if the forward orbit $\{\phi^n(x) : n\geq 0\}$
is a set of bounded (standard) height.  Thus, our main theorems
can be rephrased without reference to canonical heights by saying
that once isotrivial maps are excluded, preperiodic points in
$\PKhat$ are precisely those points whose forward iterates have
bounded height.
\end{remark}

\section{Some lemmas}
\label{sect:lem}

We will need three lemmas, two of which appeared in \cite{Ben14}.

%

\begin{lemma}
\label{lem:goodrad}
Let $\Cv$ be a complete, algebraically closed field with
non-archimedean absolute value $|\cdot|_v$.
Let $\phi\in\Cv[z]$ be a
polynomial of degree $d\geq 2$ with lead coefficient
$a_d\in\Cv^{\times}$.  Denote by $\frakK_{\phi,v}$
the filled Julia set of $\phi$ in $\Cv$, and write
$\rho_v = |a_d|_v^{-1/(d-1)}$.
Then:
\begin{list}{\rm \alph{bean}.}{\usecounter{bean}}
\item There is a unique smallest disk $U_0\subseteq\Cv$
  which contains $\frakK_{\phi,v}$.
\item $U_0$ is a closed disk of some radius
  $r_v \geq \rho_v$.
\item $\phi$ has potentially good reduction if and only
  if $r_v = \rho_v$.  In this case, $\frakK_{\phi,v}= U_0$.
\end{list}
\end{lemma}

\begin{proof}
This is the non-archimedean case of Lemma~2.5 of \cite{Ben14}.
\end{proof}


\begin{lemma}
\label{lem:manydisk}
With notation as in Lemma~\ref{lem:goodrad},
suppose that $r_v>\rho_v$.
Then $\phi^{-1}(U_0)$
is a disjoint union of closed disks $V_1,\ldots,V_{m}\subseteq U_0$,
where $2\leq m \leq d$.  Moreover, for each $i=1\ldots, m$,
$\phi(V_i)=U_0$.
\end{lemma}

\begin{proof}
This is a part of Lemma~2.7 of \cite{Ben14}.
\end{proof}

\begin{lemma}
\label{lem:smalldisks}
With notation as in Lemma~\ref{lem:goodrad},
suppose that $r_v>\rho_v$.
Let $z_0$ be any point in $\frakK_{\phi,v}$, and let $X=\phi^{-2}(z_0)$,
which is a set of at most $d^2$ points.  Then
$$\frakK_{\phi,v} \subseteq \bigcup_{x\in X} D(x,\rho_v).$$
\end{lemma}

\begin{proof}
Denote by $y_1,\ldots, y_d$ the elements (repeated according
to multiplicity) of $\phi^{-1}(z_0)$.  For each
$i=1,\ldots, d$, denote by $x_{i,1},\ldots, x_{i,d}$
the elements (again, with multiplicity)
of $\phi^{-1}(y_i)$.  Thus,
$X = \{x_{i,j} : 1\leq i,j \leq d\}$.

Let $V_1,\ldots, V_{m}$ be the disks of Lemma~\ref{lem:manydisk}.
There must be two distinct disks $V_k$ and $V_{\ell}$ of distance
exactly $r_v$ apart, by Lemma~\ref{lem:goodrad}.a and the ultrametric
triangle inequality,
or else all $m$ disks (and hence $\frakK_{\phi,v}$) would fit inside
a disk of radius smaller than $r_v$.  In addition, by
Lemma~\ref{lem:manydisk}, for each $i=1,\ldots, d$,
each of $V_k$ and $V_{\ell}$ must contain
at least one of $x_{i,1},\ldots, x_{i,d}$.  Again by the ultrametric
property, then, for every $w\in\Cv$ and every $i=1,\ldots, d$,
there is some $j=1,\ldots, d$ such that $|w-x_{i,j}|_v \geq r_v$.

Thus, given any point
$$w\in \Cv \setminus \bigcup_{x\in X} D(x,\rho_v),$$
and any index $i=1,\ldots, d$, we have
$$|\phi(w) - y_i|_v = |a_d|_v \prod_{j=1}^d |w- x_{i,j}|_v
\geq |a_d|_v \rho_v^{d-1} r_v = r_v,$$
where the inequality is by the previous paragraph, and
the final equality is by the definition of $\rho_v$.
It follows that for any such $w$,
$$|\phi^2(w) - z_0| = |a_d|_v \prod_{i=1}^d |w- y_i|_v
\geq |a_d|_v r_v^d > |a_d|_v \rho_v^{d-1} r_v = r_v.$$
However, $U_0 = \Dbar(z_0,r_v)$,
so that $\phi^2(w)\not\in U_0$, and hence $w\not\in\frakK_{\phi,v}$,
as desired.
\end{proof}

\section{Proofs}
\label{sect:pfs}

The following proposition will allow us to rephrase isotriviality
in terms of reduction.

\begin{prop}
\label{prop:iso}
Let $\FF$ and $K$ be as in Theorem~A, and let
$\phi\in K[z]$ be a polynomial of degree at least two.
Suppose that there are at least three points in $\PK$ of canonical
height zero.  Then the following are equivalent:
\begin{list}{\rm \alph{bean}.}{\usecounter{bean}}
\item There is a $K$-rational
affine change of coordinates $\gamma(z)=az+b$ for which
$\gamma^{-1}\circ\phi\circ\gamma$ is defined over $\FF_K$.
\item $\phi$ has potentially good reduction at every $v\in M_K$.
\end{list}
\end{prop}

\begin{proof}
The forward implication is immediate from the comments
following Definition~\ref{def:reduc}.
For the converse, note that $\hhat_{\phi}(\infty)=0$.
By hypothesis, there are at least two other (distinct) points
$z_0,z_1\in K$ with $\hhat_{\phi}(z_0)=\hhat_{\phi}(z_1)=0$.
Let $\gamma$ be the affine automorphism
$$\gamma(z)=(z_1-z_0)z +z_0,$$
which is defined over $K$.
Let $\psi = \gamma^{-1}\circ \phi \circ \gamma$, which is
a polynomial in $K[z]$ of the same degree as $\phi$, and
which satisfies $0,1\in \frakK_{\psi,v}$ for all $v\in M_K$.
Write $\psi(z)= b_d z^d + \cdots + b_1 z + b_0$.

By Lemma~\ref{lem:goodrad}, for each $v\in M_K$,
the filled Julia set $\frakK_{\psi,v}$
is contained in a disk of radius $\sigma_v=|b_d|_v^{-1/(d-1)}$.
In fact, since $0,1\in\frakK_{\psi,v}$,
we have $\frakK_{\psi,v}\subseteq \Dbar(0,\sigma_v)$ and $\sigma_v\geq 1$.
By the product formula
applied to $b_d$, it follows that $\sigma_v=|b_d|_v=1$ for all $v$.
Thus, by Lemma~\ref{lem:goodrad}.c, $\frakK_{\psi,v}=\Dbar(0,1)$,
and therefore,
by the comments following Definition~\ref{def:julia}, $\psi$
has good reduction at each $v\in M_K$.  By the comments following
Definition~\ref{def:reduc}, then, $|b_d|_v=1$ and $|b_i|_v\leq 1$
for all $v\in M_K$ and all $0\leq i\leq d-1$.
By the product formula, this means for each $i$ either that $b_i=0$
or that $|b_i|_v=1$ for all $v\in M_K$.  Hence, by the second property
listed near the end of Section~\ref{ssect:fields},
all coefficients of $\psi$ lie in $\FF_K$.
\end{proof}

We are now prepared to prove Theorem~A.

\begin{proof}
Let $\calZ$ denote the set of $K$-rational points points of canonical
height zero.
It suffices to show that $\calZ$ is a finite set.
Indeed, we have
already seen that all preperiodic points have canonical height zero.
Conversely, if $\calZ$ is finite, then since $\phi(\calZ)\subseteq\calZ$,
every point of $\calZ$ has finite forward orbit under $\phi$ and is
therefore preperiodic.


If $\calZ$ has fewer than three elements, we are done.
Otherwise, let $S$ be the set of places $v\in M_K$ at which $\phi$ has
bad (i.e., not potentially good) reduction, and let $s=\# S$.
As remarked in Section~\ref{sect:jul}, $s < \infty$.  By
the hypotheses and by Proposition~\ref{prop:iso}, we have $s\geq 1$.


Write $\phi(z)=a_d z^d + \cdots + a_1 z + a_0$, and for each
$v\in M_K$, let $\rho_v = |a_d|_v^{-1/(d-1)}$.
For each bad place $v\in S$, $\frakK_{\phi,v}$ is contained in
a union of at most $d^2$ open disks of radius $\rho_v$, by
Lemma~\ref{lem:smalldisks}.  In addition, 
by Lemma~\ref{lem:goodrad}, for each potentially good place
$v\in M_K\setminus S$, $\frakK_{\phi,v}$ is a
single closed disk of radius $\rho_v$.  For each choice
$\{V_v : v\in S\}$ of one of the $d^2$ disks at each bad place,
then, there is at most one $K$-rational point which lies in
$V_v$ at each $v\in S$ and in $\frakK_{\phi,v}$ at each
$v\in M_K \setminus S$.
Indeed, if there were two such points $x$ and $y$, then by the
product formula,
$$1 = \prod_{v\in S} |x-y|_v \, \cdot
\prod_{v\in M_K\setminus S} |x-y|_v
< \prod_{v \in M_K} |a_d|_v^{-1/(d-1)} = 1,$$
which is a contradiction.  (Here, the strict inequality is because
$S\neq\emptyset$.)

There are only finitely many choices of one disk $V_v$ at each $v\in S$.
Thus, there are only finitely many $K$-rational points which lie in
every filled Julia set, and hence only finitely many of canonical height
zero.
\end{proof}

\begin{remark}
\label{rem:slogs}
In the case $s\geq 1$ of Theorem~A, the above proof
shows that $\phi$ has at most $1+d^{2s}$ $K$-rational preperiodic
points, including the point at infinity.
However, we can prove a far stronger bound by invoking
verbatim the (longer and more complicated) argument in Case~1
of the proof of Theorem~7.1 in \cite{Ben14}.  That proof, which
involves bounds for certain products of differences of points
in the various filled Julia sets, shows that the number of
$K$-rational preperiodic points of $\phi$ is no more than
$$
1 + (d^2 - 2d + 2)(s\log_d s + 3s)
$$
if $1\leq s \leq d-1$, and no more than
$$ 1 + (d^2 - 2d + 2)(s\log_d s + s\log_d\log_d s + 3s)$$
if $s\geq d$.
Of course, if $s=0$ but $\phi$ still satisfies the hypotheses
of Theorem~A, then the above proof gives an upper bound of $2$.
\end{remark}

Theorem~B is now a direct consequence of Theorem~A, as follows.

\begin{proof}
If $x\in\PKhat$ is preperiodic, we already know
that $\hhat_{\phi}(x)=0$.  Conversely, if $\hhat_{\phi}(x)=0$,
let $L=K(x)$ be the field generated by $x$.
By hypothesis, since $L/K$ is a finite extension,
$\phi$ is not $L$-rationally conjugate to
a morphism defined over $\FF_L$.
Thus, by Theorem~A, all points of $\PP^1(L)$ of canonical
height zero are preperiodic.  In particular, $x$ is preperiodic.
\end{proof}

\bibliographystyle{amsplain}
\bibliography{ffheights}

\end{document}